\documentclass[reqno,centertags]{amsart}
\usepackage{amsmath}
\usepackage{amssymb}
\usepackage{amsfonts}

\theoremstyle{plain}

\newtheorem{problem}{Problem}
\newtheorem{proposition}{Proposition}

\numberwithin{equation}{section}

\input{tcilatex}

\begin{document}
\title[]{Quadric surfaces of coordinate finite type Gauss map}
\author{Hassan Al-Zoubi}
\address{Department of Basic Sciences, Al-Zaytoonah University of Jordan,
Amman}
\email{dr.hassanz@zuj.edu.jo}
\author{Tareq Hamadneh}
\address{Department of Mathematics, Al-Zaytoonah University of Jordan, P.O.
Box 130, Amman, Jordan 11733}
\email{t.hamadneh@zuj.edu.jo}
%\author{Stylianos Stamatakis}
%\address{Department of Mathematics, Aristotle University of Thessaloniki}
%\email{stamata@math.auth.gr}
\date{}

\begin{abstract}
We study quadric surfaces in the 3-dimensional Euclidean space which are of coordinate finite type Gauss map with respect to the second fundamental form $II$, i.e., their Gauss map vector $\boldsymbol{n}$ satisfies the relation $\Delta ^{II}\boldsymbol{n}=\varLambda \boldsymbol{n}$, where $\Delta ^{II}$ denotes the Laplace operator of the second fundamental form $II$ of the surface and $\varLambda$ is a square matrix of order 3. We show that helicoids and spheres are the only class of surfaces mentioned above satisfying $\Delta ^{II}\boldsymbol{n}=\varLambda \boldsymbol{n}$.
\end{abstract}

\maketitle

\section{Introduction}
In 1983 B. Y. Chen introduces the notion of Euclidean immersions of finite type, \cite{R2}, \cite{R12} and from that time on the research into surfaces of finite type has grown up as one can see in the literature in this subject. Denote by $\Delta ^{I}$ the second Beltrami-Laplace operator corresponding to the first fundamental form $I$ of a (connected) submanifold $M^{m}$ in the $n$-dimensional Euclidean space $\mathbb{E}^{n}$. Let $\boldsymbol{x},\boldsymbol{H}$ be the position vector field and the mean curvature field of $M^{m}$ respectively. Then, it is well known that \cite{R2}
\begin{equation}
\Delta ^{I}\boldsymbol{x}=-n\boldsymbol{H}.  \notag
\end{equation}

From this formula one can see that $M^{m}$ is a minimal submanifold if and only if all coordinate functions, restricted to $M^{m}$, are eigenfunctions
of $\Delta ^{I}$ with eigenvalue $\lambda =0$. Moreover in \cite{R11} T.Takahashi showed that a submanifold $M^{m}$ for which $\Delta^{I}\boldsymbol{x} =\lambda\boldsymbol{x}$, i.e. for which all coordinate functions are eigenfunctions of $\Delta ^{I}$ with the same eigenvalue $\lambda$, is either minimal in $\mathbb{E}^{n}$ with eigenvalue $\lambda=0$ or minimal in a hypersphere of $\mathbb{E}^{n}$ with eigenvalue $\lambda>0$.

The class of finite type submanifolds in an arbitrary dimensional Euclidean spaces is very large, on the other hand very little is known about surfaces of finite type in the Euclidean 3-space $\mathbb{E}^{3}$. Actually, so far, the only known surfaces of finite type corresponding to the first fundamental form in the Euclidean 3-space are the minimal surfaces, the circular cylinders and the spheres. So in \cite{R3} B.-Y. Chen laid out the problem to determine finite type surfaces in the Euclidean 3-space $\mathbb{E}^{3}$ other than the ordinary spheres and minimal ones.

In \cite{R16} O. Garay generalized T. Takahashi's condition studied surfaces in $\mathbb{E}^{3}$ for which all coordinate functions $\left(
x_{1},x_{2},x_{3}\right)$ of $\boldsymbol{x}$ satisfy $\Delta^{I}\boldsymbol{x}_{i} = \lambda_{i}x_{i}, i = 1,2,3$, not necessarily with the same
eigenvalue. Another generalization was studied in \cite{R6} for which surfaces in $\mathbb{E}^{3}$ satisfy the condition $\Delta^{I}\boldsymbol{x}= A\boldsymbol{x} + B$ $(\S)$ where $A \in\mathbb{Re}^{3\times3} ;B \in\mathbb{Re}^{3\times1}$. It was shown that a surface $S$ in $\mathbb{E}^{3}$ satisfies $(\S)$ if and only if it is an open part of a minimal surface, a sphere, or a circular cylinder. Surfaces satisfying $(\S)$ are said to be of coordinate finite type.

In this context, one can also study surfaces in $\mathbb{E}^{3}$ whose Gauss map $\boldsymbol{n}$ satisfies a relation of the form
\begin{equation}  \label{eq1.1}
\Delta^{I}\boldsymbol{n}= A\boldsymbol{n}
\end{equation}
where $A \in\mathbb{Re}^{3\times3}$.

Concerning this, in \cite{R1} it was shown that planes and circular cylinders are the only ruled surfaces whose Gauss map satisfies (\ref{eq1.1}). It is also known that \cite{D4} planes, spheres and circular cylinders are the only surfaces of revolution whose Gauss map satisfies (\ref{eq1.1}). In \cite{B6} Ch. Baikoussis and L. Verstraelen proved that planes, spheres and circular cylinders are the only helicoidal surfaces in $\mathbb{E}^{3}$ whose Gauss map satisfies (\ref{eq1.1}). Next in \cite{B5}, this problem was solved for the class of translation surfaces in $\mathbb{E}^{3}$, it was shown that the only translation surfaces whose Gauss map satisfies (\ref{eq1.1}) are the planes and the circular cylinders. Finally, in \cite{B3} this problem was studied for the compact and noncompact cyclides of Dupin, it was proved that neither for the compact nor for the noncompact cyclides of Dupin, there exist a matrix $A\in\mathbb{Re}^{3\times3}$ such that the relation (\ref{eq1.1}) is satisfied.

In this area, S. Stamatakis and H. Al-Zoubi in 2003 introduced the notion of surfaces of finite type with respect to the second or third fundamental forms \cite{S1} and it has been a topic of active research since then. As an extension of the studies mentioned above, we raise the following two questions which seem to be very interesting:

\begin{problem}
\label{2}Classifying all surfaces of coordinate finite type in the Euclidean 3-space with respect to the second or third fundamental form.
\end{problem}

\begin{problem}
\label{3}Classifying all surfaces of coordinate finite type Gauss map in the Euclidean 3-space with respect to the second or third fundamental form.
\end{problem}

%Therefore, in order to give an answer to the second and third problem, it is worthwhile investigating the classification of surfaces in the Euclidean space $\mathbb{E}^{3}$ in terms of finite $J$-type, $(J=II, III)$ by studying the families of surfaces mentioned above.
Many results concerning these two problems one can find in \cite{A1, A2, A3, A4, A5, A6, A7, A8, A9, A10, A11, A12, A13, K1, S2}.

In this paper we will pay attention to surfaces of finite $II$-type. First, we will establish a formula for $\Delta^{II}\boldsymbol{n}$ by using tensors calculations. Further, we continue our study by studying the quadric surfaces in $\mathbb{E}^{3}$ which are connected, complete and of which their Gauss map satisfy the following relation
\begin{equation}
\Delta ^{II}\boldsymbol{n}=\varLambda \boldsymbol{n},  \label{eq 4}
\end{equation}

%Moreover, in this context, the same authors proved that the surfaces $S:\boldsymbol{x}=\boldsymbol{x}(u^{1},u^{2})$ satisfying the condition
%\begin{equation*}
%\Delta ^{III}\boldsymbol{x}=\lambda \boldsymbol{x},\ \ \ \ \lambda \in\mathbb{R},
%\end{equation*}%
%i.e., for which all coordinate functions are eigenfunctions of $\Delta^{III} $ with the same eigenvalue $\lambda ,$ are precisely either the minimal surfaces $\left( \lambda =0\right) ,$ or the part of spheres $\left(\lambda =2\right) .$
%\noindent where $M(m,n)$ denotes the set of all matrices of the type $(m,n)$, while in \cite{R7} O. Garay showed that

%\begin{theorem}
%\label{TB} The only complete surfaces of revolutions in $\mathbb{R}^{3},$ whose component functions are eigenfunctions of their Laplacian, then
%the surface must be a catenoid, a sphere or a right circular cylinder.
%\end{theorem}

%Recently, H. Al-Zoubi and S. Stamatakis studied coordinate finite type with respect to the third fundamental form, more precisely, in \cite{S2} they proved the following

%\begin{theorem}
%\label{TC} A surface of revolution $S$ in $\mathbb{R}^{3},$ satisfies (\ref{eq 4}), if and only if S is a catenoid or a part of a sphere.
%\end{theorem}

\section{Preliminaries}

In the Euclidean 3-space $\mathbb{E}^{3}$ we consider a $C^{r}$-surface $S$, $r\geq3$, defined by an injective $C^{r} $-immersion $\boldsymbol{x} = \boldsymbol{x}(u^1,u^2)$ on a region $U: = I \times \mathbb{R} \,\,\,(I\subset \mathbb{R}$ open interval) of $\mathbb{R}^{2}$ whose Gaussian curvature never vanishes. We denote by
\begin{equation*}
I = g_{ij}\,d u^i d u^j, \quad II = b_{ij}\,d u^i d u^j, \quad III = e_{ij}\,d u^i d u^j, \quad i,j = 1,2,
\end{equation*}
the first, second and third fundamental forms of $S$ respectively. Let $f(u^{1}, u^{2})$ and $h(u^{1}, u^{2})$ be two sufficiently differentiable functions on $S$. Then the first differential parameter of Beltrami corresponding to the fundamental form $J = I, II, III$ is defined by \cite{H1}
\begin{equation}  \label{4.1}
\nabla^{J}(f,h):=a^{ij}f_{/i}h_{/j}
\end{equation}
where $f_{/i}: =\frac{\partial f}{\partial u^{i}}$, and $(a^{ij})$ denotes the components of the inverse tensor of $(g_{ij}), (b_{ij})$ and $(e_{ij})$ for $J = I, II$ and $III$ respectively.
The second differential parameter of Beltrami corresponding to the fundamental form $J = I, II, III$ of $S$ is defined by \cite{H1}
\begin{equation}  \label{4.2}
\triangle^{J}f:=-a^{ij}\nabla^{J}_{i}f_{/j},
\end{equation}
where $\nabla^{J}_{i}$ is the covariant derivative in the $u^{i}$ direction corresponding to the fundamental form $J$.

We compute now $\triangle^{II}\boldsymbol{n}$. Firstly, we mention the following relation for later use \cite{S1}
\begin{equation}  \label{4.21}
\nabla^{II}(h,\boldsymbol{n})+grad^{I}h=0.
\end{equation}

Applying (\ref{4.2}) for the position vector $\boldsymbol{n}$ of $S$ we have
\begin{equation}
\triangle^{II}\boldsymbol{n}=-b^{ij}(\nabla^{II}_{i})\boldsymbol{n}_{/j},
\end{equation}  \label{4.223}
Taking into consideration the equations (\cite{H1}, p.128)
\begin{equation*}
\nabla^{II}_{i}\boldsymbol{n}_{/j} =-\frac{1}{2}b^{kr}(\nabla^{III}_{r}b_{ij})\boldsymbol{n}_{/k}-e_{ij}\boldsymbol{n},
\end{equation*}
so that (\ref{4.223}) takes the form
\begin{equation*}
\triangle^{II}\boldsymbol{n}=\frac{1}{2}b^{kr}b^{ij}(\nabla^{III}_{r}b_{ij})\boldsymbol{n}_{/k}+ b^{ij}e_{ij}\boldsymbol{n}.
\end{equation*}

We consider the Christoffel symbols of the second kind corresponding to the first, second and third fundamental form, respectively
\begin{equation*}
\Gamma^{k}_{ij}:=\frac{1}{2}g^{kr}(-g_{ij/r} + g_{ir/j} + g_{jr/i}),
\end{equation*}
\begin{equation*}
\Pi^{k}_{ij}:=\frac{1}{2}b^{kr}(-b_{ij/r} + b_{ir/j} + b_{jr/i}),
\end{equation*}
\begin{equation*}
\Lambda^{k}_{ij}:=\frac{1}{2}e^{kr}(-e_{ij/r} + e_{ir/j} + e_{jr/i}),
\end{equation*}

and we put
\begin{equation}  \label{4.6}
T_{ij}^{k}: = \Gamma^{k}_{ij} -  \Pi^{k}_{ij},
\end{equation}
\begin{equation}  \label{4.61}
\widetilde{T}_{ij}^{k}: = \Lambda^{k}_{ij} -  \Pi^{k}_{ij}.
\end{equation}

It is known that (see \cite{H1}, p.22)
\begin{equation}  \label{4.7}
T_{ij}^{k}: = -\frac{1}{2}b^{kr}\nabla^{I}_{r}b_{ij},
\end{equation}
\begin{equation}  \label{4.71}
\widetilde{T}_{ij}^{k}: = -\frac{1}{2}b^{kr}\nabla^{III}_{r}b_{ij},
\end{equation}
and
\begin{equation}  \label{4.72}
\widetilde{T}_{ij}^{k}+T_{ij}^{k} = 0.
\end{equation}

On account of
\begin{equation*}
2H = b_{ij} g^{ij} = e_{ij} b^{ij},
\end{equation*}
and (\ref{4.71}) we obtain
\begin{equation*}
\triangle^{II}\boldsymbol{n}=-b^{ij}\widetilde{T}_{ij}^{k}\boldsymbol{n}_{/k}+ 2H\boldsymbol{n}.
\end{equation*}

From the Mainardi-Codazzi equations (see \cite{H1}, p.128)
\begin{equation}  \label{4.41}
\nabla^{I}_{k}b_{ij}-\nabla^{I}_{i}b_{jk}=0,
\end{equation}
and using (\ref{4.7}) and (\ref{4.72}) we have
\begin{equation}  \label{4.10}
\triangle^{II}\boldsymbol{n}=b^{kr}T_{rj}^{j}\boldsymbol{n}_{/k}+ 2H\boldsymbol{n}.
\end{equation}

For the Christoffel symbols $\Gamma^{j}_{kj}$ and $\Pi^{j}_{kj}$ we have (see \cite{H1}, p.125)
\begin{equation}  \label{4.8}
\Gamma^{j}_{ij}: = \frac{g_{/i}}{2g}, \ \ \ \ \ \,\Pi^{j}_{ij}:=\frac{b_{/i}}{2b},
\end{equation}
where $g: = det(g_{ij})$ and $b: = det(b_{ij})$.
On the other hand, the Gauss curvature $K$ of $S$ is given by
\begin{equation*}
K=\frac{b}{g}.
\end{equation*}

Hence, we have
\begin{equation}  \label{4.81}
\frac{K_{/k}}{K}=\frac{b_{/k}}{b}-\frac{g_{/k}}{g}.
\end{equation}

On the other hand using (\ref{4.6}), (\ref{4.7}), (\ref{4.8}) and (\ref{4.81}) we have
\begin{equation*}
b^{kr} T_{rj}^{j} \boldsymbol{n}_{/k}= -\frac{1}{2K}b^{kr} K_{/r} \boldsymbol{n}_{/k}=-\frac{1}{2K}\nabla^{II}(K,\boldsymbol{n}).
\end{equation*}
		
Inserting this in (\ref{4.10}) we find in view of (\ref{4.21})
\begin{equation}  \label{4.11}
\triangle^{II}\boldsymbol{n}=\frac{1}{2K}grad^{I}(K)+ 2H\boldsymbol{n}.
\end{equation}

%\section{Main result}

Now  we mention our main result which is the following

\begin{proposition}
\label{T2.2} The only quadric surfaces in the 3-dimensional \textit{Euclidean }space that satisfy (\ref{eq 4}), are the helicoids and the spheres.
\end{proposition}

Our discussion is local, which means that we show in fact that any open part of a quadric satisfies (\ref{eq 4}), if it is an open part of a helicoid or a sphere.

%Before starting the proof of our main results we first show that the surface mentioned in the above proposition indeed satisfies the condition (\ref{eq 4}).

%Let S$^{2}(r)$ be a sphere of radius $r$ centered at the origin. If $\boldsymbol{x}$ denotes the position vector field of S$^{2}(r),$ then the Gauss map $\boldsymbol{n}$\ is given by -$\frac{\boldsymbol{x}}{r}.$ For the Gauss curvature $K$ and the mean curvature $H$ of S$^{2}(r)$ we have $K$ = $\frac{1}{r^{2}}$ and $H$ = $\frac{1}{r},$ so, by virtue of (\ref{4.11}), we obtain

%\begin{equation*}
%\Delta ^{II}\boldsymbol{n}=-\frac{2}{c}\boldsymbol{n}
%\end{equation*}%
%and we find that S$^{2}(r)$ satisfies (\ref{eq 4}) with

%\begin{equation*}
%\varLambda =\left[
%\begin{array}{ccc}
%-\frac{2}{c} & 0 & 0 \\
%0 & -\frac{2}{c} & 0 \\
%0 & 0 & -\frac{2}{c}%
%\end{array}%
%\right] .
%\end{equation*}

\section{\protect\bigskip Proof of Proposition \protect\ref{T2.2}}

Let now $S$ be a quadric surface in the Euclidean 3-space $E^{3}.$ Then $S$ is either ruled, or of one of the following two kinds

\begin{equation}
z^{2}-ax^{2}-by^{2}=c,\ \ \ \ abc\neq 0  \label{I}
\end{equation}%
or

\begin{equation}
z=\frac{a}{2}x^{2}+\frac{b}{2}y^{2},\ \ \ \ a>0,\ b>0.  \label{II}
\end{equation}

We first show that helicoids are the only ruled surface in the three-dimensional Euclidean space such that its Gauss map satisfies condition (\ref{eq 4}). Next we show that the Gauss map of a quadric of the kind (\ref{I}) satisfies (\ref{eq 4}) if and only if $a=-1$ and $b=-1$, which means that $S$ is a sphere. finally, we show that the Gauss map of a quadric of the kind (\ref{II}) is never satisfying (\ref{eq 4}).

\subsection{Ruled surfaces}

Let $S$ be a ruled surface in $E^{3}.$ We suppose that $S$ is a non-cylindrical ruled surface. This surface can be expressed in terms of a directrix curve $\boldsymbol{\sigma }(s)$ and a unit vector field $\boldsymbol{\tau}(s)$ pointing along the rulings as

\begin{equation*}  \label{eq1}
S:\boldsymbol{x}(s,t)=\boldsymbol{\sigma }(s)+t\boldsymbol{\tau }(s),\ \ \
\ \ \ \ \ s\in J,\ \ \ \ \ \ \ \ -\infty <t<\infty .
\end{equation*}%
Moreover, we can take the parameter $s$ to be the arc length along the spherical curve $\boldsymbol{\tau }(s)$. Then we have

\begin{equation}  \label{eq2}
\langle \boldsymbol{\sigma }^{\prime },\boldsymbol{\tau }\rangle =0,\ \ \ \
\ \langle \boldsymbol{\tau },\boldsymbol{\tau }\rangle =1,\ \ \ \ \
\langle \boldsymbol{\tau }^{\prime },\boldsymbol{\tau }^{\prime }\rangle =1,
\end{equation}%
where the prime denotes the derivative in $s$. The first fundamental form of $S$ is

\begin{equation*}  \label{eq3}
I=qds^{2}+dt^{2},
\end{equation*}%
while the second fundamental form is

\begin{equation*}  \label{eq4}
II=\frac{p}{\sqrt{q}}ds^{2}+\frac{2A}{\sqrt{q}}dsdt,
\end{equation*}%
where

\begin{eqnarray*}  \label{eq5}
q &:&=\langle \boldsymbol{\sigma }^{\prime },\boldsymbol{\sigma }^{\prime}\rangle +2\langle \boldsymbol{\sigma }^{\prime },\boldsymbol{\tau }^{\prime }\rangle t+t^{2}, \\
p &:&=\left( \boldsymbol{\sigma }^{\prime },\boldsymbol{\tau },\boldsymbol{\sigma }^{\prime \prime }\right) +\left[ \left( \boldsymbol{\sigma }^{\prime
},\boldsymbol{\tau },\boldsymbol{\tau }^{\prime \prime }\right) +\left(\boldsymbol{\tau }^{\prime },\boldsymbol{\tau },\boldsymbol{\sigma }^{\prime \prime }\right) \right] t+\left( \boldsymbol{\tau }^{\prime },\boldsymbol{\tau },\boldsymbol{\tau }^{\prime \prime }\right) t^{2}, \\
A &:&=\left( \boldsymbol{\sigma }^{\prime },\boldsymbol{\tau },\boldsymbol{\tau }^{\prime }\right) .
\end{eqnarray*}

For convenience, we put

\begin{eqnarray*}  \label{eq6}
k&:&=\langle \boldsymbol{\sigma }^{\prime },\boldsymbol{\sigma }^{\prime }\rangle ,\ \ \ \ \ \ \ \ \ \ l :=\langle \boldsymbol{\sigma }^{\prime },\boldsymbol{\tau }^{\prime }\rangle , \\
m&:&=\left( \boldsymbol{\tau }^{\prime },\boldsymbol{\tau },\boldsymbol{\tau }^{\prime \prime }\right) ,\ \ \ \ \ n :=\left( \boldsymbol{\sigma }%
^{\prime },\boldsymbol{\tau },\boldsymbol{\tau }^{\prime \prime }\right)+\left( \boldsymbol{\tau }^{\prime },\boldsymbol{\tau },\boldsymbol{\sigma
}^{\prime \prime }\right) , \\
r&:&=\left( \boldsymbol{\sigma }^{\prime },\boldsymbol{\tau },\boldsymbol{\sigma }^{\prime \prime }\right) ,
\end{eqnarray*}
and thus we have

\begin{equation*}  \label{eq7}
q=t^{2}+2lt+k ,\ \ \ \ \ \ \ \ \ \ p=m t^{2}+n t+r .
\end{equation*}

Furthermore, the Gaussian curvature $K$ of $S$ is given by

\begin{equation}  \label{eq8}
K=-\frac{A^{2}}{q^{2}}.
\end{equation}

Since the Gaussian curvature $K$ never vanishes. Thus we have
\begin{equation}     \label{eq9}
A\neq 0.
\end{equation}

A parametric representation of the Gauss map $\boldsymbol{n}$ of $S$ is given by
\begin{equation}  \label{eq10}
\boldsymbol{n }(s,t) =\frac{1}{\sqrt{q}}\big(\boldsymbol{\sigma}'\times\boldsymbol{\tau} + t\boldsymbol{\tau}'\times\boldsymbol{\tau}\big).
\end{equation}

If, for simplicity, we put
\begin{equation*}
\boldsymbol{P}:=\boldsymbol{P}(s)=\boldsymbol{\sigma}'\times\boldsymbol{\tau}, \ \ \ \ \ \boldsymbol{Q}:=\boldsymbol{Q}(s)=\boldsymbol{\tau}'\times\boldsymbol{\tau},
\end{equation*}
then the Gauss map $\boldsymbol{n}$ of $S$ becomes
\begin{equation}  \label{eq11}
\boldsymbol{n }(s,t) =\frac{1}{\sqrt{q}}\big(\boldsymbol{P} + t\boldsymbol{Q}\big).
\end{equation}

The Beltrami operator with respect to the second fundamental form, after a lengthy computation, can be expressed as follows

\begin{eqnarray}  \label{eq12}
\Delta ^{II} =\frac{2\sqrt{q}}{A}\frac{\partial ^{2}}{\partial s\partial t}-\frac{p\sqrt{q}}{A^{2}} \frac{\partial ^{2}}{\partial t^{2}}
-\frac{\sqrt{q}p_{t}}{A^{2}}\frac{\partial }{\partial t},
\end{eqnarray}%
where
\begin{equation*}
p_{t}:=\frac{\partial p}{\partial t}
\end{equation*}

Applying (\ref{eq12}) for the position vector $\boldsymbol{n}$ one finds

\begin{eqnarray}  \label{eq13}
\Delta ^{II}\boldsymbol{n} &=&\frac{1}{q^{2}}\bigg[\bigg(\frac{tpqq_{tt}}{2A^{2}}-\frac{3tpq_{t}^{2}}{4A^{2}}+\frac{tqq_{t}p_{t}}{A^{2}}-\frac{q^{2}p_{t}}{A^{2}}+\frac{qpq_{t}}{A^{2}}\\
&&+\frac{3tq_{t}q_{s}}{2A}-\frac{tqq_{st}}{A}-\frac{qq_{s}}{A}\bigg)\boldsymbol{Q}+\bigg(\frac{2q^{2}}{A}-\frac{tqq_{t}}{A}\bigg)\boldsymbol{Q}'-  \frac{qq_{t}}{A}\boldsymbol{P}'  \notag   \\
&&+\bigg(\frac{pqq_{tt}}{2A^{2}}-\frac{3pq_{t}^{2}}{4A^{2}}+\frac{qq_{t}p_{t}}{2A^{2}}+\frac{3q_{s}q_{t}}{2A}-\frac{qq_{st}}{A}\bigg)\boldsymbol{P}\bigg]. \notag
\end{eqnarray}
Here again we have
\begin{equation*}
p_{s}:=\frac{\partial p}{\partial s}, \ \ \ q_{t}:=\frac{\partial q}{\partial t},  \ \ \ q_{s}:=\frac{\partial q}{\partial s}
\end{equation*}
and the prime stands, as we mention before, the derivative with respect to $s$, that is
\begin{equation*}
\boldsymbol{P}'=\frac{d\boldsymbol{P}}{d s}, \ \ \ \boldsymbol{Q}'=\frac{d\boldsymbol{Q} }{d s}.
\end{equation*}

Equation (\ref{eq13}) can be expressed as follows
\begin{eqnarray}  \label{eq14}
\Delta ^{II}\boldsymbol{n}&=&\frac{1}{q^{2}}\bigg[\frac{1}{A^{2}}f_{1}(t)\boldsymbol{Q(s)}+\frac{1}{A}f_{2}(t)\boldsymbol{Q'(s)} \\ &&+\frac{1}{A^{2}}f_{3}(t)\boldsymbol{P(s)}+\frac{1}{A}f_{4}(t)\boldsymbol{P'(s)}\bigg] \notag
\end{eqnarray}
where $f_{i}(t), i= 1,2,3,4$ are polynomials in $t$ with functions in $s$ as coefficients, and deg$(f_{i})\leq4$. More precisely, we have

\begin{eqnarray} \label{eq15}
f_{1}(t) &=&2mt^{5}+\big(n+4lm\big)t^{4} \\
&&+ \big(3kl-3l^{2}m+2ln+2l'A\big)t^{3}  \notag \\
&&+\big(3kn-3l^{2}n-4klm+2k'A-2ll'A\big)t^{2} \notag \\
&&+\big(3kr-3l^{2}r-2kln-2k^{2}m+k'lA-4kl'A\big)t \notag \\
&&-k^{2}n-kk'A, \notag
\end{eqnarray}
\begin{eqnarray} \label{eq16}
f_{2}(t) =2lt^{3}+\big(4l^{2}n+2k\big)t^{2} + 6klt+2k^{2},
\end{eqnarray}
\begin{eqnarray} \label{eq17}
f_{3}(t) &=&\big(2lm-n\big)t^{3} + \big(4l'+l^{2}m-ln-2r+3km\big)t^{2} \\
&&+\big(2kn-l^{2}n+2klm+3k'A3+2ll'-4lr\big)t \notag \\
&&+kr+kln-2kl'-3l^{2}r+3k'l, \notag
\end{eqnarray}
\begin{eqnarray} \label{eq18}
f_{4}(t) =2t^{3}+ 6lt^{2}+\big(4l^{2}+2k\big)t+2kl.
\end{eqnarray}

Let now $\boldsymbol{n}=(n_{1},n_{2},n_{3}),\boldsymbol{P }=(P_{1},P _{2},P _{3})$ and $\boldsymbol{Q }=(Q _{1},Q_{2},Q _{3})$ be the coordinate functions of $\boldsymbol{n},\boldsymbol{P}$ and $\boldsymbol{Q}$ respectively. By virtue of (\ref{eq14}) we obtain
\begin{eqnarray} \label{eq19}
\Delta ^{II}n_{i}&=&\frac{1}{q^{2}}\bigg[\frac{1}{A^{2}}f_{1}(t)Q_{i}(s)+\frac{1}{A}f_{2}(t)Q'_{i}(s) \\ &&+\frac{1}{A^{2}}f_{3}(t)P_{i}(s)+\frac{1}{A}f_{4}(t)P'_{i}(s)\bigg], \\ \notag
i &=&1,2,3.  \notag
\end{eqnarray}

We denote by $\lambda _{ij},i,j=1,2,3$ the entries of the matrix $\varLambda $. Using (\ref{eq19}) and condition (\ref{eq 4}) we have

\begin{eqnarray*}
&&\frac{1}{A^{2}}f_{1}(t)Q_{i}(s)+\frac{1}{A}f_{2}(t)Q'_{i}(s)+\frac{1}{A^{2}}f_{3}(t)P_{i}(s)+\frac{1}{A}f_{4}(t)P'_{i}(s) \\
&=&\lambda _{i1}q^{\frac{3}{2}}\big(P_{i} + tQ_{i}\big)+\lambda _{i2}q^{\frac{3}{2}}\big(P_{i} + tQ_{i}\big)+\lambda _{i3}q^{\frac{3}{2}}\big(P_{i} + tQ_{i}\big), \\
i &=&1,2,3.
\end{eqnarray*}

Consequently%
\begin{eqnarray}
&&-3m ^{2}\tau _{i}t^{5}+[\left( m ^{\prime }A-m A^{\prime }-4m n
-7\lambda m ^{2}\right) \tau _{i}+3m A\tau _{i}^{\prime }]t^{4}  \notag
\\
&&+[m A\sigma _{i}^{\prime }-A^{2}\tau _{i}^{\prime \prime }+(2n
A+7\lambda m A)\tau _{i}^{\prime }+(n ^{\prime }A-n A^{\prime }
\notag \\
&&+2\lambda m ^{\prime }A-2\lambda m A^{\prime }-\lambda ^{\prime }m
A-A^{2}-10\lambda m n -2m r -n ^{2}  \notag \\
&&-4k m ^{2})\tau _{i}]t^{3}+[(k m ^{\prime }A-k m
A^{\prime }-\frac{1}{2}k ^{\prime }m A+2\lambda n ^{\prime }A
\notag \\
&&-2\lambda n A^{\prime }-\lambda ^{\prime }n A-r A^{\prime }+r
^{\prime }A-3\lambda A^{2}-3\lambda n ^{2}-6\lambda m r   \notag \\
&&-6k m n )\tau _{i}+3\lambda m A\sigma _{i}^{\prime }-2\lambda
A^{2}\tau _{i}^{\prime \prime }-A^{2}\sigma _{i}^{\prime \prime }+(\lambda
^{\prime }A+5\lambda n   \notag \\
&&+4k m +r )A\tau _{i}^{\prime }]t^{2}+[(\frac{1}{2}k
^{\prime }A+3k n +3\lambda r )A\tau _{i}^{\prime }  \notag \\
&&+(k n ^{\prime }A-k n A^{\prime }-\frac{1}{2}k ^{\prime
}n A+2\lambda r ^{\prime }A-2\lambda r A^{\prime }-\lambda ^{\prime
}r A  \notag \\
&&-k A^{2}-2\lambda ^{2}A^{2}-2k n ^{2}+r ^{2}-2\lambda n
r -4k m r )\tau _{i}  \notag \\
&&-2\lambda A^{2}\sigma _{i}^{\prime \prime }-k A^{2}\tau _{i}^{\prime
\prime }+\left( \lambda n -r +2k m +\lambda ^{\prime }A\right)
A\sigma _{i}^{\prime }]t  \notag \\
&&-A^{4}\left( \lambda _{i1}\tau _{1}+\lambda _{i2}\tau _{2}+\lambda
_{i3}\tau _{3}\right) t+(k r ^{\prime }A-k r A^{\prime }
\notag \\
&&-\frac{1}{2}k ^{\prime }r A+\lambda r ^{2}-k \lambda
A^{2}-2k n r )\tau _{i}-k A^{2}\sigma _{i}^{\prime \prime }
\notag \\
&&+2k r A\tau _{i}^{\prime }+(\frac{1}{2}k ^{\prime }A-\lambda
r +k n )A\sigma _{i}^{\prime }  \notag \\
&&-A^{4}\left( \lambda _{i1}\sigma _{1}+\lambda _{i2}\sigma _{2}+\lambda
_{i3}\sigma _{3}\right) =0.  \label{eq 10}
\end{eqnarray}

For $i=1,2,3,$ (\ref{eq 10}) is a polynomial in $t$ with functions in $s$ as
coefficients. This implies that the coefficients of the powers of $t$ in (%
\ref{eq 10}) must be zeros, so we obtain, for $i=1,2,3,$ the following
equations

\begin{equation}
3m ^{2}\tau _{i}=0,  \label{eq 11}
\end{equation}

\begin{equation*}
\left( m ^{\prime }A-m A^{\prime }-4m n -7\lambda m ^{2}\right)
\tau _{i}+3m A\tau _{i}^{\prime }=0,
\end{equation*}%
\begin{eqnarray}
&&m A\sigma _{i}^{\prime }-A^{2}\tau _{i}^{\prime \prime }+(2n
A+7\lambda m A)\tau _{i}^{\prime }  \notag \\
&&+(n ^{\prime }A-n A^{\prime }+2\lambda m ^{\prime }A-2\lambda m
A^{\prime }-\lambda ^{\prime }m A  \notag \\
&&-A^{2}-10\lambda m n -2m r -n ^{2}-4k m ^{2})\tau
_{i}=0,  \label{eq 12}
\end{eqnarray}%
\begin{eqnarray}
&&(k m ^{\prime }A-k m A^{\prime }-\frac{1}{2}k ^{\prime
}m A+2\lambda n ^{\prime }A-2\lambda n A^{\prime }-\lambda ^{\prime
}n A  \notag \\
&&-r A^{\prime }+r ^{\prime }A-3\lambda A^{2}-3\lambda n
^{2}-6\lambda m r -6k m n )\tau _{i}+3\lambda m A\sigma
_{i}^{\prime }  \notag \\
&&-2\lambda A^{2}\tau _{i}^{\prime \prime }-A^{2}\sigma _{i}^{\prime \prime
}+(\lambda ^{\prime }A+5\lambda n +4k m +r )A\tau _{i}^{\prime
}=0,  \label{eq 13}
\end{eqnarray}%
\begin{eqnarray}
&&(k n ^{\prime }A-k n A^{\prime }-\frac{1}{2}k ^{\prime
}n A+2\lambda r ^{\prime }A-2\lambda r A^{\prime }-\lambda ^{\prime
}r A-k A^{2}  \notag \\
&&-2\lambda ^{2}A^{2}-2k n ^{2}+r ^{2}-2\lambda n r -4k
m r )\tau _{i}-2\lambda A^{2}\sigma _{i}^{\prime \prime }-k
A^{2}\tau _{i}^{\prime \prime }  \notag \\
&&+(\frac{1}{2}k ^{\prime }A+2k n +4\lambda r )A\tau
_{i}^{\prime }+\left( \lambda n -r +2k m +\lambda ^{\prime
}A\right) A\sigma _{i}^{\prime }  \notag \\
&=&A^{4}(\lambda _{i1}\tau _{1}+\lambda _{i2}\tau _{2}+\lambda _{i3}\tau
_{3}),  \label{eq 14}
\end{eqnarray}%
\begin{eqnarray}
&&(k r ^{\prime }A-k r A^{\prime }-\frac{1}{2}k
^{\prime }r A+\lambda r ^{2}-k \lambda A^{2}-2k n r
)\tau _{i}+2k r A\tau _{i}^{\prime }  \notag \\
&&+(\frac{1}{2}k ^{\prime }A-\lambda r +k n )A\sigma
_{i}^{\prime }-k A^{2}\sigma _{i}^{\prime \prime }  \notag \\
&=&A^{4}\left( \lambda _{i1}\sigma _{1}+\lambda _{i2}\sigma _{2}+\lambda
_{i3}\sigma _{3}\right) .  \label{eq 15}
\end{eqnarray}

From (\ref{eq 11}) one finds

\begin{equation}
m =\left( \boldsymbol{\tau }^{\prime },\boldsymbol{\tau },\boldsymbol{%
\tau }^{\prime \prime }\right) =0,  \label{eq 16}
\end{equation}%
which implies that the vectors $\boldsymbol{\tau }^{\prime },\boldsymbol{%
\tau },\boldsymbol{\tau }^{\prime \prime }$ are linearly dependent, and
hence there exist two functions $\sigma _{1}=\sigma _{1}(s)$ and $\sigma
_{2}=\sigma _{2}(s)$ such that

\begin{equation}
\boldsymbol{\tau }^{\prime \prime }=\sigma _{1}\boldsymbol{\tau +}\sigma
_{2}\boldsymbol{\tau }^{\prime }.  \label{eq 17}
\end{equation}

On differentiating $\langle \boldsymbol{\tau }^{\prime },\boldsymbol{\tau }%
^{\prime }\rangle =1,$ we obtain $\langle \boldsymbol{\tau }^{\prime },%
\boldsymbol{\tau }^{\prime \prime }\rangle =0.$ So from (\ref{eq 17}) we
have

\begin{equation}
\boldsymbol{\tau }^{\prime \prime }=\sigma _{1}\boldsymbol{\tau }.
\label{eq 19}
\end{equation}

By taking the derivative of $\langle \boldsymbol{\tau },\boldsymbol{\tau }%
\rangle =1$ twice, we find that

\begin{equation*}
\langle \boldsymbol{\tau }^{\prime },\boldsymbol{\tau }^{\prime }\rangle
+\langle \boldsymbol{\tau },\boldsymbol{\tau }^{\prime \prime }\rangle =0.
\end{equation*}

But $\langle \boldsymbol{\tau }^{\prime },\boldsymbol{\tau }^{\prime
}\rangle =1,$ and taking into account (\ref{eq 19}) we find that $\sigma
_{1}(s)=-1.$ Thus (\ref{eq 19}) becomes $\boldsymbol{\tau }^{\prime \prime
}=-\boldsymbol{\tau }$ which implies that

\begin{equation}
\tau _{i}^{\prime \prime }=-\tau _{i},\ \ \ \ i=1,2,3.  \label{eq 20}
\end{equation}

Using (\ref{eq 16})\ and (\ref{eq 20})\ equation (\ref{eq 12})\ reduces to

\begin{equation*}
2n A\tau _{i}^{\prime }+(n ^{\prime }A-n A^{\prime }-n ^{2})\tau
_{i}=0,\ \ \ \ i=1,2,3
\end{equation*}%
or, in vector notation

\begin{equation}
2n A\boldsymbol{\tau }^{\prime }+(n ^{\prime }A-n A^{\prime }-n
^{2})\boldsymbol{\tau }=\boldsymbol{0}.  \label{eq 21}
\end{equation}

By taking the derivative of $\langle \boldsymbol{\tau },\boldsymbol{\tau }%
\rangle =1$, we find that the vectors $\boldsymbol{\tau },\boldsymbol{\tau
}^{\prime }$ are linearly independent, and so from (\ref{eq 21}) we obtain
that $n A=0.$ We note that $A\neq 0$, since from (\ref{eq8}) the Gauss
curvature vanishes, so we are left with $n =0$. Then equation (\ref{eq 13}%
) becomes

\begin{equation*}
-A^{2}\sigma _{i}^{\prime \prime }+(\lambda ^{\prime }A+r )A\tau
_{i}^{\prime }+(r ^{\prime }A-r A^{\prime }-\lambda A^{2})\tau
_{i}=0,\ \ \ \ i=1,2,3
\end{equation*}%
or, in vector notation

\begin{equation}
-A^{2}\boldsymbol{\sigma }^{\prime \prime }+(\lambda ^{\prime }A+r )A%
\boldsymbol{\tau }^{\prime }+(r ^{\prime }A-r A^{\prime }-\lambda
A^{2})\boldsymbol{\tau }=\boldsymbol{0}.  \label{eq 22}
\end{equation}

Taking the inner product of both sides of the above equation with $%
\boldsymbol{\tau }^{\prime }$ we find in view of (\ref{eq2}) that

\begin{equation}
-A^{2}\langle \boldsymbol{\sigma }^{\prime \prime },\boldsymbol{\tau }%
^{\prime }\rangle +r A+\lambda ^{\prime }A^{2}=0.  \label{eq 23}
\end{equation}

On differentiating $\lambda =\langle \boldsymbol{\sigma }^{\prime },%
\boldsymbol{\tau }^{\prime }\rangle $ with respect to $s,$ by virtue of (%
\ref{eq 20}) and (\ref{eq2}), we get

\begin{equation}
\lambda ^{\prime }=\langle \boldsymbol{\sigma }^{\prime \prime },\boldsymbol{%
\tau }^{\prime }\rangle +\langle \boldsymbol{\sigma }^{\prime },\boldsymbol{%
\tau }^{\prime \prime }\rangle =\langle \boldsymbol{\sigma }^{\prime \prime
},\boldsymbol{\tau }^{\prime }\rangle -\langle \boldsymbol{\sigma }^{\prime
},\boldsymbol{\tau }\rangle =\langle \boldsymbol{\sigma }^{\prime \prime },%
\boldsymbol{\tau }^{\prime }\rangle .  \label{eq 24}
\end{equation}

Hence, (\ref{eq 23}) reduces to $r A=0,$ which implies that $r =0.$
Thus the vectors $\boldsymbol{\sigma }^{\prime },\boldsymbol{\tau },%
\boldsymbol{\sigma }^{\prime \prime }$ are linearly dependent, and so there
exist two functions $\sigma _{3}=\sigma _{3}(s)$ and $\sigma _{4}=\sigma
_{4}(s)$ such that

\begin{equation}
\boldsymbol{\sigma }^{\prime \prime }=\sigma _{3}\boldsymbol{\tau +}\sigma
_{4}\boldsymbol{\sigma }^{\prime }.  \label{eq 25}
\end{equation}

Taking the inner product of both sides of the last equation with $%
\boldsymbol{\tau }$ we find in view of (\ref{eq2}) that $\sigma
_{3}=\langle \boldsymbol{\sigma }^{\prime \prime },\boldsymbol{\tau }%
\rangle .$

Now, by taking the derivative of $\langle \boldsymbol{\sigma }^{\prime },%
\boldsymbol{\tau }\rangle =0,$ we find $\langle \boldsymbol{\sigma }%
^{\prime \prime },\boldsymbol{\tau }\rangle +\langle \boldsymbol{\sigma }%
^{\prime },\boldsymbol{\tau }^{\prime }\rangle =0,$ that is

\begin{equation}
\langle \boldsymbol{\sigma }^{\prime \prime },\boldsymbol{\tau }\rangle
+\lambda =0,  \label{eq 26}
\end{equation}%
and hence $\sigma _{3}=-\lambda .$

Taking again the inner product of both sides of equation (\ref{eq 25}) with $%
\boldsymbol{\tau }^{\prime }$ we find in view of (\ref{eq2}) that

\begin{equation}
\langle \boldsymbol{\sigma }^{\prime \prime },\boldsymbol{\tau }^{\prime
}\rangle =\sigma _{4}\lambda .  \label{eq 27}
\end{equation}

Using (\ref{eq 24}) we find $\lambda ^{\prime }=\sigma _{4}\lambda .$ Thus $%
\sigma _{4}=\frac{\lambda ^{\prime }}{\lambda }.$ Therefore

\begin{equation}
\boldsymbol{\sigma }^{\prime \prime }=-\lambda \boldsymbol{\tau +}\frac{%
\lambda ^{\prime }}{\lambda }\boldsymbol{\sigma }^{\prime }.  \label{eq 28}
\end{equation}

We distinguish two cases

\textbf{Case 1}: $\lambda =0.$ Because of $r =0$ equation (\ref{eq 22})
would yield $A=0,$ which is clearly impossible for the surfaces under
consideration.

\textbf{Case 2}: $\lambda \neq 0.$ From (\ref{eq 22}), (\ref{eq 28}) and $%
r =0$ we find that

\begin{equation*}
-\frac{\lambda ^{\prime }}{\lambda }A^{2}\boldsymbol{\sigma }^{\prime
}+\lambda ^{\prime }A^{2}\boldsymbol{\tau }^{\prime }=\mathbf{0}
\end{equation*}%
which implies that $\lambda ^{\prime }(\boldsymbol{\sigma }^{\prime
}-\lambda \boldsymbol{\tau }^{\prime })=\mathbf{0}.$

If $\lambda ^{\prime }\neq 0,$ then $\boldsymbol{\sigma }^{\prime }=\lambda
\boldsymbol{\tau }^{\prime }.$ Hence $\boldsymbol{\sigma }^{\prime },%
\boldsymbol{\tau }^{\prime }$ are linearly dependent, and so $A=0$ which
contradicts our previous assumption. Thus $\lambda ^{\prime }=0.$ From (\ref%
{eq 28}) we have

\begin{equation}
\boldsymbol{\sigma }^{\prime \prime }=-\lambda \boldsymbol{\tau }.
\label{eq 29}
\end{equation}

On the other hand, by taking the derivative of $k $ and using the last
equation we obtain that $k $ is constant. Hence equations (\ref{eq 14})
and (\ref{eq 15}) reduce to

\begin{eqnarray*}
\lambda _{i1}\tau _{1}+\lambda _{i2}\tau _{2}+\lambda _{i3}\tau _{3} &=&0,
\\
\lambda _{i1}\sigma _{1}+\lambda _{i2}\sigma _{2}+\lambda _{i3}\sigma _{3}
&=&0,\ \ i=1,2,3\
\end{eqnarray*}%
and so $\lambda _{ij}=0,i,j=1,2,3.$

Since the parameter $s$ is the arc length of the spherical curve $%
\boldsymbol{\tau }(s)$, and because of (\ref{eq 16}) we suppose, without
loss of generality, that the parametrization of $\boldsymbol{\tau }(s)$ is

\begin{equation*}
\boldsymbol{\tau }(s)=(\cos s,\sin s,0).
\end{equation*}

Integrating (\ref{eq 29}) twice we get

\begin{equation*}
\boldsymbol{\sigma }(s)=(c_{1}s+c_{2}+\lambda \cos s,c_{3}s+c_{4}+\lambda
\sin s,c_{5}s+c_{6}),
\end{equation*}%
where $c_{i},i=1,2,...,6$ are constants.

Since $k =\langle \boldsymbol{\sigma }^{\prime },\boldsymbol{\sigma }%
^{\prime }\rangle $ is constant, it's easy to show that $c_{1}=c_{3}=0.$
Hence $\boldsymbol{\sigma }(s)$ reduces to

\begin{equation*}
\boldsymbol{\sigma }(s)=(c_{2}+\lambda \cos s,c_{4}+l \sin
s,c_{5}s+c_{6}).
\end{equation*}

Thus we have

\begin{equation*}
S:\boldsymbol{x}(s,t)=(c_{2}+(l +t)\cos s,c_{4}+(l +t)\sin
s,c_{5}s+c_{6})
\end{equation*}%
which is a helicoid.

\subsection{Quadrics of the first kind}

\noindent This kind of quadric surfaces can be parametrized as follows

\begin{equation*}
\boldsymbol{x}(u,v)=\left( u,v,\sqrt{c+au^{2}+bv^{2}}\right) .
\end{equation*}

Let's denote the function $c+au^{2}+bv^{2}$ by $\omega $. Then, using the natural frame $\{{\boldsymbol{x}_{u}, \boldsymbol{x}_{v}}\}$ of $S$ defined by
\begin{equation*}
\boldsymbol{x_{u}}=\left( 1,0,\frac{au}{\sqrt{\omega}}\right),
\end{equation*}
and
\begin{equation*}
\boldsymbol{x_{v}}=\left( 0,1,\frac{bv}{\sqrt{\omega}} \right),
\end{equation*}
the components $g_{ij}$ of the first fundamental form in (local) coordinates are the following

\begin{equation*}
g_{11}=<\boldsymbol{x_{u}},\boldsymbol{x_{u}}>= 1+\frac{\left( au\right) ^{2}}{\omega },\ \ \
\end{equation*}
\begin{equation*}
g_{12}=<\boldsymbol{x_{u}},\boldsymbol{x_{v}}>=\frac{abuv}{\omega },\ \ \
\end{equation*}
\begin{equation*}
g_{22}=<\boldsymbol{x_{v}},\boldsymbol{x_{v}}>=1+\frac{\left( bv\right) ^{2}}{\omega },
\end{equation*}

Denoting by $\boldsymbol{n}$ the Gauss map of $S$, then we have
\begin{equation*}
\boldsymbol{n}=\Big({-\frac{au}{\sqrt{\Phi}},-\frac{bv}{\sqrt{\Phi}},\frac{\sqrt{\omega}}{\sqrt{\Phi}}}\Big),
\end{equation*}
where $\Phi= c+a(a+1)u^{2}+b(b+1)v^{2}$. The components $b_{ij}$ of the second fundamental form in (local) coordinates are the following
\begin{equation*}
b_{11}=<\boldsymbol{x_{uu}},\boldsymbol{n}>= \frac{a\left( bv^{2}+c\right)}{\omega\sqrt{\Phi} },\ \ \
\end{equation*}
\begin{equation*}
b_{12}=<\boldsymbol{x_{uv}},\boldsymbol{n}>=-\frac{abuv}{\omega\sqrt{\Phi} },\ \ \
\end{equation*}
\begin{equation*}
b_{22}=<\boldsymbol{x_{vv}},\boldsymbol{n}>=\frac{b\left( au^{2}+c\right)}{\omega\sqrt{\Phi} },
\end{equation*}

By using the natural frame $\{{\boldsymbol{n}_{u}, \boldsymbol{n}_{v}}\}$ of $S$ defined by
\begin{equation*}
\boldsymbol{n_{u}}= \big( -\frac{a[b(b+1)v^{2}+c]}{\Phi^{\frac{3}{2}}},\frac{ab(a+1)uv}{\Phi^{\frac{3}{2}}},\frac{au[b(b+1)v^{2}-ac]}{\sqrt{\omega}\Phi^{\frac{3}{2}}}\big) ,
\end{equation*}
and
\begin{equation*}
\boldsymbol{n_{v}}=\big( \frac{ab(b+1)uv}{\Phi^{\frac{3}{2}}},-\frac{b[a(a+1)u^{2}+c]}{\Phi^{\frac{3}{2}}},\frac{bv[a(a+1)u^{2}-bc]}{\sqrt{\omega}\Phi^{\frac{3}{2}}}\big).
\end{equation*}

The second Beltrami differential operator with respect to the second fundamental form is defined by
\begin{equation*}
\triangle^{II}f = \frac{-1}{\sqrt{|b|}} \, \frac{\partial \left(\sqrt{|b|} \,b^{ij} \frac{\partial f}{\partial u^i}\right)}{\partial u^j},
\end{equation*}
where $b : =\det(b_{ij})$. After a long computation the Beltrami operator $\Delta ^{II}$ of $S$ can be expressed as follows:
\begin{eqnarray}  \label{4}
\triangle ^{II} &=&-\frac{\sqrt{\Phi}}{c}\Bigg[\frac{au^{2}+c}{a}\frac{\partial ^{2}} {\partial u^{2}} - 2uv \frac{\partial ^{2}}{\partial u\partial v} +\frac{bv^{2}+c}{b}{\partial v^{2}}  \notag \\
&&+2u\frac{\partial }{\partial u}+2v\frac{\partial }{\partial v}\Bigg],
\end{eqnarray}

Applying (\ref{4}) for the functions $n_{1}= -\frac{au}{\sqrt{\Phi}}$ and $n_{2}= -\frac{bv}{\sqrt{\Phi}}$, one finds

\begin{equation}  \label{eq 32}
\Delta ^{II}n_{1}=\Delta ^{II}\big(-\frac{au}{\sqrt{\Phi}}\big)= \frac{au}{\Phi^{2}}\big[3b(b+1)(b-a)v^{2}-3ac-(b+2)\Phi \big].
\end{equation}
\begin{equation}  \label{eq 33}
\Delta ^{II}n_{2}=\Delta ^{II}\big(-\frac{bv}{\sqrt{\Phi}}\big)= \frac{bv}{\Phi^{2}}\big[3a(a+1)(a-b)u^{2}-3bc-(a+2)\Phi \big].
\end{equation}

We denote by $\lambda _{ij},i,j=1,2,3$ the entries of the matrix $\Lambda $. On account of (\ref{4}) we get

\begin{equation}  \label{35}
\Delta ^{II}n_{1}=\Delta ^{II}\Big(-\frac{au}{\sqrt{\Phi}}\Big)=\lambda _{11}\Big(-\frac{au}{\sqrt{\Phi}}\Big)+\lambda _{12}\Big(-\frac{bv}{\sqrt{\Phi}}\Big)+\lambda _{13}\Big(\frac{\sqrt{\omega}}{\sqrt{\Phi}}\Big),
\end{equation}

\begin{equation}  \label{36}
\Delta ^{II}n_{2}=\Delta ^{II}\Big(-\frac{bv}{\sqrt{\Phi}}\Big)=\lambda _{21}\Big(-\frac{au}{\sqrt{\Phi}}\Big)+\lambda _{22}\Big(-\frac{bv}{\sqrt{\Phi}}\Big)+\lambda _{23}\Big(\frac{\sqrt{\omega}}{\sqrt{\Phi}}\Big),
\end{equation}

We denote by $l _{ij},i,j=1,2,3$ the entries of the matrix $\varLambda $.
On account of (\ref{eq 4}) we get

\begin{equation}
\Delta ^{III}x_{1}=\Delta ^{III}u=l _{11}u+l _{12}v+l _{13}%
\sqrt{\omega },  \label{eq 32}
\end{equation}

\begin{equation}
\Delta ^{III}x_{2}=\Delta ^{III}v=l _{21}u+l _{22}v+l _{23}%
\sqrt{\omega },  \label{eq 33}
\end{equation}

\begin{equation*}
\Delta ^{III}x_{3}=\Delta ^{III}\sqrt{\omega }=l _{31}u+l
_{32}v+l _{33}\sqrt{\omega }.
\end{equation*}

Applying (\ref{eq 31}) on the coordinate functions $x_{i},i=1,2$ of the
position vector $\boldsymbol{x}$ and by virtue of (\ref{eq 32}) and (\ref{eq
33}), we find respectively

\begin{eqnarray}
\Delta ^{III}u &=&-\frac{uT}{c^{2}}\left[ 3\left( a+1\right) u^{2}+3\left(
b+1\right) v^{2}+\frac{c(3b+a+2ab)}{ab}\right]  \notag \\
&=&l _{11}u+l _{12}v+l _{13}\sqrt{\omega },  \label{eq 34}
\end{eqnarray}

\begin{eqnarray}
\Delta ^{III}v &=&-\frac{vT}{c^{2}}\left[ 3\left( a+1\right) u^{2}+3\left(
b+1\right) v^{2}+\frac{c(b+3a+2ab)}{ab}\right]  \notag \\
&=&l _{21}u+l _{22}v+l _{23}\sqrt{\omega },  \label{eq 35}
\end{eqnarray}

Putting $v=0$ in (\ref{eq 34}), we obtain that%
\begin{eqnarray*}
&&-\frac{3a(a+1)^{2}}{c^{2}}u^{5}-\frac{(a+1)(6b+a+2ab)}{bc}u^{3}-\frac{%
\left( 3b+a+2ab\right) }{ab}u \\
&=&l _{11}u+l _{13}\sqrt{c+au^{2}}.
\end{eqnarray*}

Since $a\neq 0$ and $c\neq 0$ this implies that $a=-1.$

Similarly, if we put $u=0$ in (\ref{eq 35}) we obtain that%
\begin{eqnarray*}
&&-\frac{3b(b+1)^{2}}{c^{2}}v^{5}-\frac{(b+1)(b+6a+2ab)}{ac}v^{3}-\frac{%
\left( b+3a+2ab\right) }{ab}v \\
&=&l _{22}v+l _{23}\sqrt{c+bv^{2}}.
\end{eqnarray*}

This implies that $b=-1.$ Hence $S$ must be a sphere.

\subsection{Quadrics of the second kind}

For this kind of surfaces we can consider a parametrization

\begin{equation*}
\boldsymbol{x}(u,v)=\left( u,v,\frac{a}{2}u^{2}+\frac{b}{2}v^{2}\right) .
\end{equation*}

Then the components $g_{ij},b_{ij}$ and $e_{ij}$ of the first, second and
third fundamental tensors are the following

\begin{equation*}
g_{11}=1+\left( au\right) ^{2},\ \ \ g_{12}=abuv,\ \ \ g_{22}=1+\left(
bv\right) ^{2},
\end{equation*}

\begin{equation*}
b_{11}=\frac{a}{\sqrt{g}},\ \ \ b_{12}=0,\ \ \ b_{22}=\frac{b}{\sqrt{g}},
\end{equation*}

\begin{equation*}
e_{11}=\frac{a^{2}}{g^{2}}(1+b^{2}v^{2}),\ \ \ e_{12}=-\frac{a^{2}b^{2}}{%
g^{2}}uv,\ \ \ e_{22}=\frac{b^{2}}{g^{2}}(1+a^{2}u^{2}),
\end{equation*}%
where $g:=\det \left( g_{ij}\right) =1+\left( au\right) ^{2}+\left(
bv\right) ^{2}.$

A straightforward computation shows that the Beltrami operator $\Delta
^{III} $ of $S$ takes the following form

\begin{eqnarray}
\Delta ^{III} &=&-\frac{g(1+a^{2}u^{2})}{a^{2}}\frac{\partial ^{2}}{\partial
u^{2}}-\frac{g(1+b^{2}v^{2})}{b^{2}}\frac{\partial ^{2}}{\partial v^{2}}
\notag \\
&&-2uvg\frac{\partial ^{2}}{\partial u\partial v}-2ug\frac{\partial }{%
\partial u}-2vg\frac{\partial }{\partial v}.  \label{eq 36}
\end{eqnarray}

On account of (\ref{eq 4}) we get

\begin{equation}
\Delta ^{III}x_{1}=\Delta ^{III}u=l _{11}u+l _{12}v+l
_{13}\left( \frac{a}{2}u^{2}+\frac{b}{2}v^{2}\right) ,  \label{eq 37}
\end{equation}

\begin{equation}
\Delta ^{III}x_{2}=\Delta ^{III}v=l _{21}u+l _{22}v+l
_{23}\left( \frac{a}{2}u^{2}+\frac{b}{2}v^{2}\right) ,  \label{eq 38}
\end{equation}

\begin{equation*}
\Delta ^{III}x_{3}=\Delta ^{III}\sqrt{\omega }=l _{31}u+l
_{32}v+l _{33}\left( \frac{a}{2}u^{2}+\frac{b}{2}v^{2}\right) .
\end{equation*}

Applying (\ref{eq 36}) on the coordinate functions $x_{i},i=1,2$ of the
position vector $\boldsymbol{x}$ and by virtue of (\ref{eq 37}) and (\ref{eq
38}) we find respectively

\begin{equation}
\Delta ^{III}u=-2ug=l _{11}u+l _{12}v+l _{13}\left( \frac{a%
}{2}u^{2}+\frac{b}{2}v^{2}\right) ,  \label{eq 39}
\end{equation}

\begin{equation}
\Delta ^{III}v=-2vg=l _{21}u+l _{22}v+l _{23}\left( \frac{a%
}{2}u^{2}+\frac{b}{2}v^{2}\right) .  \label{eq 40}
\end{equation}

Putting $v=0$ in (\ref{eq 39}), we obtain that

\begin{equation*}
-2a^{2}u^{3}-2u=l _{11}u+l _{13}\frac{a}{2}u^{2}.
\end{equation*}

This implies that $a$ must be zero$.\ $Putting $u=0$ in (\ref{eq 40}), we
obtain that
\begin{equation*}
-2b^{2}v^{3}-2v=l _{22}v+l _{23}\frac{b}{2}v^{2}.
\end{equation*}
This implies that $b$ must be zero$,$ which is clearly impossible, since $%
a>0 $ and$\ b>0.$

%\noindent \textbf{Acknowledgement}\\
%\noindent The authors would like to express their thanks to the referee for his useful remarks.

\end{document}